\documentclass{amsart}

\usepackage[latin1]{inputenc}
\usepackage{color}
\usepackage{dsfont}
\usepackage{varioref}
\usepackage{epsfig}
\usepackage{hyperref}
\usepackage[english]{babel}

\newtheorem{Lem}{Lemma}
\newtheorem{Obs}{Observation}
\newtheorem{Thm}{Theorem}
\newtheorem{Cor}[Thm]{Corollary}

\newcommand{\dP}{{\mathds P}}
\newcommand{\dC}{{\mathds C}}
\newcommand{\dR}{{\mathds R}}
\newcommand{\dQ}{{\mathds Q}}
\newcommand{\dF}{{\mathds F}}
\newcommand{\smcdot}{{\textup{$\cdot$}}}

\begin{document}

\title{A Septic with 99 real Nodes}
\author{Oliver Labs}
\address{Johannes Gutenberg Universit\"at Mainz,
  Germany}
\email{Labs@Mathematik.Uni-Mainz.de, mail@OliverLabs.net}
\date{September 9, 2004}


\begin{abstract}
We find a surface of degree $7$ in $\dP^3(\dR)$ with
$99$ real nodes within a family of surfaces with dihedral symmetry: 
First, we consider this family over some small prime fields, 
which allows us to test all possible parameter sets using computer algebra.
In this way we find some examples of $99$-nodal surfaces over some of these 
finite fields.
Then, the examination of the geometry of these surfaces allows us to
determine the parameters of a $99$-nodal septic in characteristic zero. 
This narrows the possibilities for $\mu(7)$, the maximum number of nodes on a
septic, to: $99 \le \mu(7) \le 104$. 
When reducing our surface modulo $5$, we even obtain a $100$-nodal septic
in $\dP^3(\dF_5)$.
\end{abstract}

\maketitle

\section*{Introduction}

We denote by $\mu(d)$ the maximum number of nodes (i.e.\ singularities of type
$A_1$, also called ordinary double points) a surface of degree $d$ in
$\dP^3(\dC)$ can have.  
The restrictions on $\mu(d)$ known so far are summarized in the following
table:  
\begin{center}
\begin{tabular}{|r|c|c|c|c|c|c|c|c|c|c|c|c|}
\hline
&&&&&&&&&&&&\\[-1em]
degree & 2 & 3 & 4 & 5 & 6 & 7 & 8 & 9 & 10 & 11 & 12 & d\\[0.1em] 
\hline
\hline
&&&&&&&&&&&&\\[-0.9em]
$\mu(d) \ge$ & 1 & 4 & 16 & 31 & 65 & 93 & 168 & 216 & 345 & 425 & 
600 & ${5 \over 12}d^3$\\[0.25em] 
\hline
&&&&&&&&&&&&\\[-0.9em]
$\mu(d) \le$ & 1 & 4 & 16 & 31 & 65 & 104 & 174 & 246 & 360 & 480
& 645 & ${4\over9}d^3$\\[0.25em] 
\hline
\end{tabular}
\end{center}
In this article we show:
\begin{equation}\label{mu7}\mu(7) \ge 99.\end{equation}

From the table above, we see in particular that $\mu(d)$ is known up to
$d=6$. 
The upper bound $\mu(7) \le 104$ is given by Varchenko's spectrum bound 
\cite{varBound}. 
Note that for $d=7$ Miyaoka's bound \cite{miyP3} is $112$, 
but Givental's bound \cite{givBound} also gives $104$. 

The previously known septic with the greatest number of nodes was the example 
of Chmutov \cite{chmuP3} with $93$ nodes: It comes
from a construction that works for any degree $d$.  
For $d\le 5$ and the even degrees $d = 6, 8, 10, 12$ there are examples
exceeding Chmutov's lower bound: \cite{bar65}, \cite{endrOct},
\cite{ales12}.   
These had been obtained by using some beautiful geometric arguments based on
Segre's idea \cite{BSegCon1}. 

In this note, we explain how to use computer algebra experiments over prime 
fields to treat the case $d=7$ and to find the first surface of odd degree
greater than $5$ that exceeds Chmutov's general lower bound.  
Given an explicit equation of a family of hypersurfaces, there is
in fact an algorithm to find those examples with the greatest number of nodes: 
We already applied this succesfully in \cite{labsSext35}, but because of
computer performance restrictions we cannot use this 
technique in the present case. 
Instead, we choose a more geometric approach to study the family.

I thank D.~van~Straten for his permanent motivation and many valuable
discussions. 
Furthermore, I thank W.~Barth for his invitation to Erlangen which was a good
motivation to complete this work. 
I thank Slawomir Cynk for helpful discussions.
Finally, I thank St.~Endra\ss{} for discussions, motivation and
his Ph.D.\ thesis which is a great source for dihedral-symmetric surfaces
with many singularities.

\section{The Family}
\label{secTheFamily}

Inspired by many authors (see in particular: \cite{bar65},
\cite{endrThesis}, \cite{endrOct}), we look for septics with many nodes in
$\dP^3(\dC)$ within a $7$-parameter family of surfaces $S_{a_1, a_2, \dots,
  a_7} := P - U_{a_1, a_2, \dots, a_7}$ of degree $7$ admitting the dihedral
symmetry $D_7$ of a $7$-gon:
\begin{eqnarray*}P &:=& 2^6 \cdot \Pi_{j=0}^{6}\left[\cos\left({2\pi j\over
        7}\right)x 
    + \sin\left({2\pi j\over 7}\right)y -z\right]\\
& = & x \smcdot \left[x^{6}-3\smcdot7\smcdot x^{4}y^{2}
+5\smcdot7\smcdot x^{2}y^{4}-7\smcdot y^{6}\right]\\
& & + 7 \smcdot z \smcdot\!\left[\left(x^2+y^2\right)^3
 - 2^3\smcdot z^{2} \smcdot\!\left(x^2+y^2\right)^2
+ 2^4\smcdot z^{4} \smcdot\!\left(x^{2}+y^{2}\right)\right]
- 2^6 \smcdot z^{7},\\
  U_{a_1, a_2, \dots, a_7} &:=& (z+a_5
  w)\left(a_1z^3+a_2z^2w+a_3zw^2+a_4w^3+(a_6z+a_7w)(x^2+y^2)\right)^2.
\end{eqnarray*}

$P$ is the product of $7$ planes in $\dP^3(\dC)$ meeting in the point
$(0:0:0:1)$ and admitting $D_7$-symmetry with rotation axes $\{x=y=0\}$: 
In fact, $P$ is invariant under the map $y \mapsto -y$ and $P \cap 
\{z=z_0\}$ is a regular $7$-gon for $z_0 \ne 0$. 
$U$ is also $D_7$-symmetric, because $x$ and $y$ only appear as $x^2+y^2$.  

The generic surface $S$ has nodes at the $3\cdot21=63$
intersections of the ${7\choose2}=21$ doubled lines of $P$ with the doubled
cubic of $U$. 
We look for parameters $a_1, a_2, \dots, a_7$, s.t.\ the corresponding
surface has $99$ nodes.

As $S_{a_1, a_2, \dots, a_7}(x,y,z,\lambda w) = 
S_{a_1, \lambda a_2, \lambda^2 a_3, \lambda^3 a_4, \lambda a_5, a_6, \lambda
  a_7}(x,y,z,w) \ \forall \lambda \in 
\dC^*$, we choose $a_7 := 1$.
Moreover, experiments over prime fields suggest that the maximum number of
nodes on such surfaces is $99$ and that such examples exist for $a_6=1$. 
As we are mainly interested in finding an example with $99$ nodes, we restrict
ourselves to the sub-family:   
$$S := S_{a_1, a_2, a_3, a_4, a_5, 1, 1} = P - U_{a_1, a_2, a_3, a_4, a_5, 1,
  1}.$$ 
Some other cases, e.g.\ $a_6=0$, also lead to $99$-nodal septics; this 
will be discussed elsewhere.

\section{Reduction to the Case of Plane Curves}

To simplify the problem of locating examples with $99$ nodes within
our family $S$, we restrict our attention to
the $\{y=0\}$ -plane and search for plane curves $S|_{y=0}$ (we write $S_y$
for short) with many nodes.
This is possible because of the symmetry of the construction
(see \cite[p.~18, cor.~2.3.10]{endrThesis} for details):
\begin{Lem}\label{lemDdsym}
A member $S = S_{a_1, a_2, a_3, a_4, a_5, 1, 1}$ of our family of surfaces has 
only ordinary double points as singularities, if $(1:i:0:0) \notin S$ and the
surface does only contain ordinary double points as singularities in the plane
$\{y=0\}$.   
If the plane septic $S_y$ has exactly $n$ nodes and if exactly $n_{xy}$
of these nodes are on the axes $\{x=y=0\}$ then the surface $S$ has exactly
$n_{xy} + 7\cdot(n-n_{xy})$ nodes and no other singularities.
Each singularity of $S_y$ which is not on $\{x=y=0\}$ gives an orbit of $7$
singularities of $S$ under the action of the dihedral group $D_7$. 
\end{Lem}

To understand the geometry of the plane septic $S_y$ better, we 
look at the singularities that occur for generic values of the parameters.
First, we compute:
\begin{eqnarray*}
P|_{y=0} &=& x^7+7\cdot x^6z-7\cdot2^3\cdot x^4z^3+7\cdot2^4\cdot
x^2z^5-2^6\cdot z^7\\ 
&=& {(x-z)\over2^4} \cdot \big(\underbrace{x+(-\rho)z}_{=: L_1}\big)^2 \cdot
\big(\underbrace{2x+(\rho^2+4\rho)z}_{=: L_2}\big)^2 \cdot 
\big(\underbrace{2x+(-\rho^2-2\rho+8)z}_{=: L_3}\big)^2,\\ 
U|_{y=0} &=& (z+a_5
w){\big(\underbrace{(z+w)x^2+a_1z^3+a_2z^2w+a_3zw^2+a_4w^3}_{=:
    C}\big)}^2,
\end{eqnarray*}
where $\rho$ satisfies: 
\begin{equation}\label{eqnrho}\rho^3 + 2^2\rho^2-2^2\rho-2^3=0.\end{equation} 

\begin{figure}[htbp]
\begin{center}
\framebox[1.05\width]{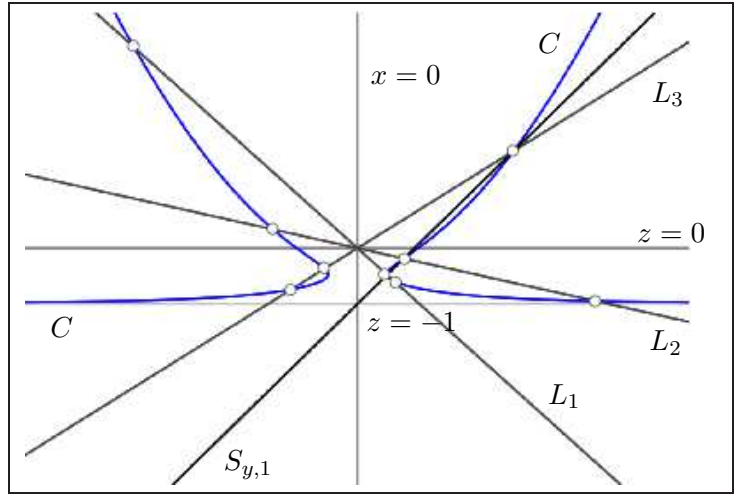}
\end{center}
\caption{The three doubled lines $L_i$ and the doubled cubic $C$
  intersect in $3\cdot3=9$ points $G_{ij}$. These are the \emph{generic}
  singularities of the plane septic $S_y$.}
\label{fig_dLdC}
\end{figure}

The three points $G_{ij}$ of intersection of $C$ with the line $L_i$
are ordinary double points of the plane septic $S_y = P|_{y=0}-U|_{y=0}$
for generic values of the parameters, s.t.\ we have $3\cdot3=9$ \emph{generic} 
singularities (see fig.~\vref{fig_dLdC}). 

\section{Finding Solutions over some Prime Fields}

By running over all possible parameter combinations over some
small prime fields $\dF_p$ using the computer algebra system {\sc
  Singular} \cite{Singular}, we find some $99$-nodal surfaces over these
fields: 
For a given set of parameters $a_1, a_2, \dots, a_5$, we can easily check the
actual number of nodes on the corresponding surface using
computer algebra (see \cite[appendix A, p.~487]{gpSingComAlg}).   

As indicated in the previous section, we work in the plane $\{y=0\}$ for
faster computations. 
It turns out that the greatest number of nodes on $S_y$ is $15$ over the
small prime fields $\dF_{p}$, \ $11 \le p \le 53$: See table \vref{tabym1}.
The prime fields $\dF_p$, $2 \le p \le 7$, are not listed because they are
special cases: These primes appear as coefficients or exponents in the
equation of our family.   
In each of the cases we checked, one of the $15$ singular points lies on the  
axes $\{x=0\}$, such that the corresponding surface has exactly
$14\cdot7+1=99$ nodes and no other singularities.

\section{The Geometry of the $15$-nodal septic Plane Curve}  

To find parameters $a_1, a_2, \dots, a_5$ in characteristic $0$ we want to use 
geometric properties of the $15$-nodal septic plane curve $S_y$.  
But as we do not know any such property yet, we use our prime field examples
to get some good ideas:
\begin{Obs}\label{obsrefined}
In all our prime field examples of $15$-nodal plane septics $S_y$, we have: 
\begin{enumerate}
\item
$S_y$ splits into a line $S_{y,1}$ and a sextic $S_{y,6}$: \ $S_y = S_{y,1}
\cdot S_{y,6}$.  
The plane curve $S_{y,6}$ of degree $6$ has $15-6=9$
singularities. Note that this property is similar to the one of the $31$-nodal
  $D_5$-symmetric quintic in $\dP^3(\dC)$ constructed by W.~Barth: See  
  \cite[p.~27-32]{endrThesis} for a description.
\end{enumerate}
The line and the sextic have some interesting geometric properties (see
fig.\ \vref{fig_dLdC} and fig.\ \vref{fig15nodes}): 
\begin{enumerate}
\setcounter{enumi}{1}
\item $S_{y,1} \, \cap \, S_{y,6} = \{R, G_{1j_1}, G_{2j_2}, G_{3j_3}, O_1,
  O_2\}$, 
  where $R$ is a point on the axes $\{x=0\}$ and the 
  $G_{ij_k}$ are three of the $9$ generic singularities $G_{ij}$ of $S_y$, one
  on each line $L_i$, and $O_1, O_2$ are some other points that neither lie on
  $\{x=0\}$, nor on one of the $L_i$.
\item The sextic $S_{y,6}$ has the six generic 
  singularities $G_{ij}, \ (i,j)\in\{1,2,3\}^2 \ \backslash$
  $\{(1,j_1),(2,j_2),(3,j_3)\}$, and three 
  exceptional singularities: $E_1, E_2, E_3$. 
\end{enumerate}
In many prime field experiments, we have furthermore:
\begin{enumerate}
\setcounter{enumi}{3}
\item In the projective $x,z,w$-plane, the point $R$ has the coordinates
  $(0:-1:1)$, s.t.\ the line $S_{y,1}$ has the form $S_{y,1}: \  z + t \cdot x
  + w = 0$ for some parameter $t$ (see also table
\vref{tabym1}).  
\end{enumerate}
The other cases ($R = (0:c:1), c\ne-1$) lead to more complicated equations and
will not be discussed here.
\end{Obs}
  
\begin{table}[htbp]
\begin{center}
\begin{tabular}{|r|c|c|c|c|c|c|c|l|}
\hline
Field & $a_1$ & $a_2$ & $a_3$ & $a_4$ & $a_5$ & $S_{y,1}$ & $\alpha$\\  
\hline
\hline
\rule{0pt}{1.1em}$\dF_{11}$ & 2 & 3 & 5 & 2 & -5 & $z=x-w$ & $\alpha=-3$ \\[0.2em]
\hline
\rule{0pt}{1.1em}$\dF_{19}$ & -7 & -2 & 7 & 1 & 8 & $z=8x-w$ & $\alpha=7$\\
$\dF_{19}$ & 2 & 0 & 1 & 9 & 7 & $z=9x-w$ & $\alpha=-4$\\
$\dF_{19}$ & 5 & -9 & 7 & -3 & -1 & $z=2x-w$ & $\alpha=-3$\\[0.2em]
\hline
\rule{0pt}{1.1em}$\dF_{23}$ & -5 & 11 & 10 & 1 & 7 & $z=-9x-w$ & $\alpha=-2$\\[0.2em]
\hline
\rule{0pt}{1.1em}$\dF_{31}$ & -15 & -13 & -5 & 13 & -10 & $z=-2x-w$ & $\alpha=-13$\\
$\dF_{31}$ & 1 & -2 & 14 & -9 & 11 & $z=15x-w$ & $\alpha=-11$\\
$\dF_{31}$ & 14 & -10 & -13 & -14 & -11 & $z=-13x-w$ & $\alpha=-7$\\[0.2em]
\hline
\rule{0pt}{1.1em}$\dF_{43}$ & -11 & 15 & 0 & -13 & -6 & $z=-6x-w$ & $\alpha=7$\\
$\dF_{43}$ & 20 & 16 & -1 & -14 & 10 & $z=-12x-w$ & $\alpha=14$\\
$\dF_{43}$ & -9 & 3 & -3 & -11 & 5 & $z=18x-w$ & $\alpha=-21$\\[0.2em]
\hline
\rule{0pt}{1.1em}$\dF_{53}$ & -8 & 20 & 14 & 18 & 11 & $z=25x-w$ & $\alpha=4$\\
$\dF_{53}$ & -2 & -10 & -14 & -26 & 16 & $z=-9x-w$ & $\alpha=24$\\
$\dF_{53}$ & 10 & 25 & -4 & 22 & 25 & $z=-16x-w$ & $\alpha=25$\\[0.2em]
\hline
\end{tabular}
\caption{A few examples of parameters giving $15$-nodal septic plane curves
(and $99$-nodal surfaces) over prime fields (see \cite{labsAlgSurf} for more
exhaustive tables).} 
\label{tabym1}
\end{center}
\end{table}

Using this observation as a guess for our septic in characteristic
$0$, we obtain several polynomial conditions on the parameters. 
Using {\sc Singular} to eliminate variables, we find the following
relation between the parameters $a_4$ and $t$:   
\begin{equation}\label{eqn_a4t}t \cdot \big(\underbrace{a_4t^3+t}_{=:
      \alpha}\big)^2 + t-1=0,\end{equation}  
which can be parametrized by $\alpha$: $t = -{1\over 1+\alpha^2}, \ a_4 =
(\alpha(1+\alpha^2)-1)(1+\alpha^2)^2.$ 
Further eliminations allow us to express all the other parameters in terms of  
$\alpha$:   
\begin{itemize}
\item
$a_1 = \alpha^7+7\alpha^5-\alpha^4+7\alpha^3-2\alpha^2-7\alpha-1,$
\item
$a_2 = (\alpha^2+1)(3\alpha^5+14\alpha^3-3\alpha^2+7\alpha-3),$
\item
$a_3 = (\alpha^2+1)^2(3\alpha^3+7\alpha-3),$
\item
$a_5 = -{\alpha^2 \over 1+\alpha^2}.$
\end{itemize}

\section{The $1$-parameter Family of Plane Sextics}

We use once more our explicit examples of $15$-nodal
septic plane curves over prime fields to finally be able to write down a
condition for $\alpha$ in characteristic $0$.

First, note that we can now easily obtain the equation of $S_{y,6}$ by
dividing the equation of our septic curve $S_y$ by the equation of the line 
$S_{y,1} = z+tx+w = z - {1\over1+\alpha^2}x +w$. 
$S_{y,6}$ is a sextic which has $6$ nodes for generic $\alpha$, but
should have $9$ double points for some special values of $\alpha$. 
One idea to determine these particular values is to find a geometric
relation between the $6$ generic singular points and the $3$ exceptional ones. 

\subsection{Three Conics}
Looking at the equations describing the singular points of our examples of
$9$-nodal sextics $S_{y,6}$ over the prime fields, we see the following:
\begin{Obs}\label{obs3conics}
For all our $9$-nodal examples of plane sextics over prime fields, there
are three conics through six of these points each 
(see fig.~\vref{fig3conics}): 
\begin{figure}[htbp]
\begin{center}
\framebox[1.05\width]{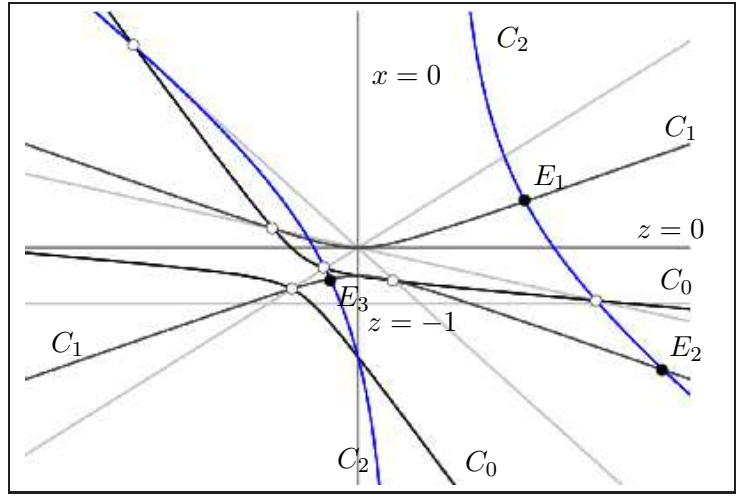}
\caption{Three conics relating the $9$ double points of the sextic
  $S_{y,6}$. $E_1, E_2,$ and $E_3$ (black) are the exceptional singularities
  (i.e.\ they do not lie on one of the lines $L_i$, see
  fig.~\vref{fig_dLdC}).  
The white points are the generic singularities, coming from the intersection
  of the doubled cubic $C$ with the three doubled lines $L_i$.}  
\label{fig3conics}
\end{center}
\end{figure}
\begin{enumerate}
\item one conic $C_0$ through the $6$ generic singularities,
\item one conic $C_1$ through the $3$ exceptional singularities and $3$ of the
  generic ones, 
\item one conic $C_2$ through the $3$ exceptional singularities and the other $3$
  generic ones.
\end{enumerate} 
Moreover, the three conics have the following
properties over the prime fields: 
\begin{enumerate}
\setcounter{enumi}{3}
\item $C_1$ has the form: 
\begin{equation}\label{eqnC1}C_1: \ x^2+kz^2+(k+4)zw=0,\end{equation}
where $k$ is a still unknown parameter. 
In particular, $C_1$ is symmetric with respect to $x\mapsto -x$ and contains
the point $(0:0:1)$. 
\item $C_0$ intersects the other two conics on the $\{x=0\}$ -axes (see
  fig.~\vref{fig3conics}):  
\begin{equation}\label{eqnConicsx0} X_1 := C_0 \cap C_1 \cap \{x=0\}, 
\quad X_2 := C_0 \cap C_2 \cap \{x=0\}.\end{equation}  
\end{enumerate}
\end{Obs}

To determine the new parameter $k$ in equation (\ref{eqnC1}), we
will use (\ref{eqnConicsx0}). 
We compute the two points of
$C_0$ on the $\{x=0\}$ -axes explicitly using {\sc Singular}:
First, the ideal $I_{S_{y,6}}^{gen}$ describing the six generic singularities
of $S_{y,6}$ can be computed from the ideal $I_{S_y}^{gen} := (C, \,
L_1L_2L_3)$  
describing the $9$ generic singularities of $S_y$ by calculating the following 
ideal quotient: $I_{S_{y,6}}^{gen} = I_{S_y}^{gen} : S_{y,1}$. 
Now, the equation of $C_0$ can be obtained by taking the degree-$2$-part of
the ideal $I_{S_{y,6}}^{gen}$:  
\begin{equation}\label{eqnC0}
\begin{array}{rrcl}C_0: \ & 
\alpha x^{2}+(\alpha ^{3}+5\alpha -1)xz+(\alpha ^{5}+(\alpha ^{3}+\alpha
-1)xw\\
&+6\alpha ^{3}-\alpha ^{2}+\alpha -1)z^{2}+(2\alpha ^{5}+8\alpha ^{3}-2\alpha
^{2}+6\alpha -2)zw\\
&+(\alpha ^{5}+2\alpha ^{3}-\alpha ^{2}+\alpha
-1)w^2&=&0.\end{array}\end{equation}   

Thus, $\{P^+, P^-\} := C_0 \cap \{x=0\} = \left\{\left(0 :
  {-2(\alpha^3+3\alpha-1)(1+\alpha^2) 
  \pm \beta(\alpha)  
\over 2(\alpha^5+6\alpha^3-\alpha^2+\alpha-1)} : 1\right)\right\},$ where 
\begin{equation}\label{eqnbeta}
\begin{array}{rcl}\beta(\alpha)^2 &:=& (\alpha^3+3\alpha-1)^2(1+\alpha^2)^2\\
&& - 4(\alpha^5+6\alpha^3-\alpha^2+\alpha-1)
  (1+\alpha^2)(\alpha^3+\alpha-1)\\[0.2em]
&=& 16\alpha(2\alpha^5+4\alpha^3-\alpha^2+2\alpha-1).
\end{array}\end{equation}  

$C_1$ intersects the $\{x=0\}$ -axes in exactly two points: $(0:0:1)$ and
$X_1$.
Hence, we can determine the two possibilities for the
parameter $k\in\dQ(\alpha,\beta(\alpha))$ in equation (\ref{eqnC1}) for $C_1$:
Together with the $z$ and $w$-coordinates of the points $P^\pm$, $C_1 \cap
\{x=0\} = \{kz^2 + kzw +4zw = 0\}$ leads to the following two possibilities:   
\begin{equation}
C_1: \quad x^2 + {-4 P_z^\pm \over P_z^\pm (P_z^\pm + 1)} z(z+w) +
4zw = 0.\end{equation}

\subsection{The Condition on $\alpha$}

The equations of the conics $C_0$ and $C_1$ will allow us to compute the condition on $\alpha$, s.t.\ the sextic $S_{y,6}$
has $9$ singularities, using the following (see observation \ref{obs3conics}
and fig.~\ref{fig3conics}): 
\begin{itemize}
\item
$C_0$ intersects the three doubled lines $L_i$ exactly in the six generic
singularities. 
\item
$C_1$ intersects the three doubled lines $L_i$ exactly in three of these six
generic singularities and the origin (which counts three times).
\end{itemize}

Thus, the set of $z$-coordinates of the three points $(C_1
\cap L_1L_2L_3) \ \backslash \ \{(0:0:1)\}$ has to be contained in the
set of $z$-coordinates of the six points $C_0 \cap L_1L_2L_3$. 
This means that the remainder $q$ of the following division ($res_x$ denotes 
the resultant with respect to $x$)
\begin{equation}\label{eqnresx}res_x(C_0, \, L_1L_2L_3)  =
  p(z) \cdot \left({1\over z^3} \cdot res_x(C_1, \, L_1L_2L_3)\right) +
  q(z)\end{equation}  
should vanish: $q = 0$.

As the degree of the remainder is $\deg(q) = 2$, this gives $3$ conditions on
$\alpha$ and $\beta(\alpha)$, coming from the fact that all the $3$ coefficients of
$q(z)$ have to vanish. 
It turns out that it suffices to take one of these, the coefficient of
$z^2$, which can be written in the form $c(\alpha)+\beta(\alpha) d(\alpha)$,
where 
$c(\alpha)$ and $d(\alpha)$ are polynomials in $\dQ[\alpha]$.  
As a condition on $\alpha$ only we can take: $$cond(\alpha) := 
\big(c(\alpha)+\beta(\alpha)d(\alpha)\big) \cdot \big(c(\alpha)-\beta(\alpha)
d(\alpha)\big) \in \dQ[\alpha],$$ 
which is of degree $150$.

This condition $cond(\alpha)$ vanishes for those $\alpha$ for which the
corresponding surface has $99$ nodes and for several other $\alpha$. 
To obtain a condition which exactly describes those $\alpha$ we are looking  
for, we factorize $cond(\alpha) = f_1 \cdot f_2 \cdots f_k$ (e.g., 
using {\sc Singular} again). 
Substituting in each of these factors our solutions over the prime fields, we
see that the only factor that vanishes is: $7\alpha^3+7\alpha+1=0.$
\begin{figure}[htbp]
\begin{center}
\framebox[1.05\width]{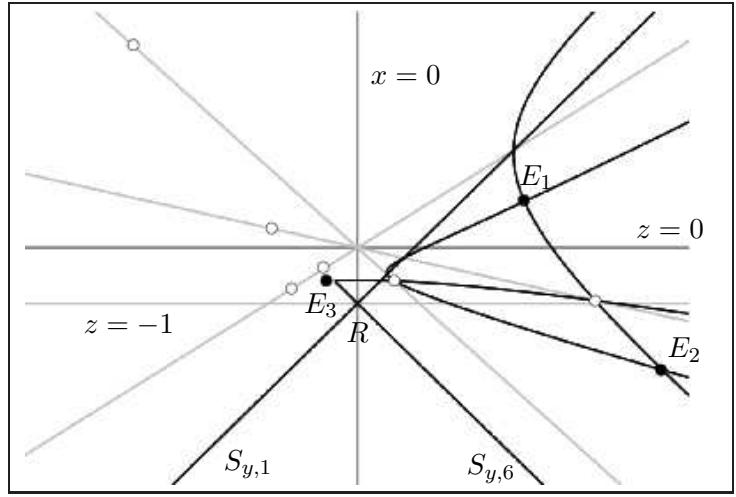}
\caption{The $15$-nodal plane septic $S_{y_{\alpha_\dR}} =
  S_{{y,1}_{\alpha_\dR}} \cdot S_{{y,6}_{\alpha_\dR}}$ (see
  (\ref{eqnalphaR}) on page \pageref{eqnalphaR}); the singularities of the
  sextic 
  $S_{{y,6}_{\alpha_\dR}}$ are marked by large circles: The three exceptional
  singularities $E_1, E_2, E_3$ are marked in black, the generic singularities
  in white.  
  The five left-most nodes are real isolated ones. 
  Only five of the six intersections of the line $S_{{y,1}_{\alpha_\dR}}$ and
  the sextic $S_{{y,6}_{\alpha_\dR}}$ are visible because we just show a small
  part of the whole $(x,z)$-plane.}   
\label{fig15nodes}
\end{center}
\end{figure}

\section{The Equation of the $99$-nodal Septic}

Up to this point, it is still only a guess --- verified over some prime fields
--- that the values $\alpha$ satisfying the condition above give
$99$-nodal septics in characteristic $0$.
But a straightforward computation with {\sc Singular} shows: 
\begin{Thm}[$99$-nodal Septic]
Let $\alpha\in\dC$ satisfy:
\begin{equation}\label{eqnalpha}7\alpha^3+7\alpha+1=0.\end{equation}
Then the surface $S_\alpha$ in $\dP^3(\dC)$ of degree $7$ with equation 
$S_\alpha := P-U_\alpha$ has exactly $99$ ordinary double points and no other
singularities, where  
\begin{eqnarray*}P &:=& x \smcdot \left[x^{6}-3\smcdot7\smcdot x^{4}y^{2}
+5\smcdot7\smcdot x^{2}y^{4}-7\smcdot y^{6}\right]\\
& & + 7 \smcdot z \smcdot\!\left[\left(x^2+y^2\right)^3
 - 2^3\smcdot z^{2} \smcdot\!\left(x^2+y^2\right)^2
+ 2^4\smcdot z^{4} \smcdot\!\left(x^{2}+y^{2}\right)\right]
- 2^6 \smcdot z^{7},\\
  U_\alpha &:=& (z+a_5
  w)\left((z+w)(x^2+y^2)+a_1z^3+a_2z^2w+a_3zw^2+a_4w^3\right)^2, 
\end{eqnarray*}
\begin{center}\begin{tabular}{l@{\qquad}l}$a_1 :=
  -{12\over7}\alpha^2-{384\over49}\alpha-{8\over7}$, 
  & $a_2 := -{32\over7}\alpha^2+{24\over49}\alpha-4,$\\[0.2em] 
$a_3 := -4\alpha^2+{24\over49}\alpha-4$, &  
$a_4 := -{8\over7}\alpha^2+{8\over49}\alpha-{8\over7},$\\[0.2em] 
$a_5 := 49\alpha^2-7\alpha+50$.
\end{tabular}\end{center}
There is exactly one real solution $\alpha_\dR \in \dR$ to the condition 
(\ref{eqnalpha}), 
\begin{equation}\label{eqnalphaR}\alpha_\dR \approx 
  -0.14010685,\end{equation} 
and all the singularities of $S_{\alpha_\dR}$ are also real.  
\end{Thm}
\begin{proof}
To show that the surface has $99$ nodes and no other singularities 
we verify the hypothesis of lemma \vref{lemDdsym}, 
because this speeds up the computations. 
The algorithm is straightforward.
Our {\sc Singular} code is listed in the appendix and can also be obtained
from \cite{labsAlgSurf}. 

Using the fact that $\beta(\alpha)^2 = \left({12\over7}\right)^2$ together
with the geometric description of the singularities of the plane septic given
in the previous sections, it is also straightforward to verify the reality
assertion.  
\end{proof}

%
\section{Concluding Remarks}

\begin{figure}[htbp]
\begin{center}
\includegraphics[width=3in]{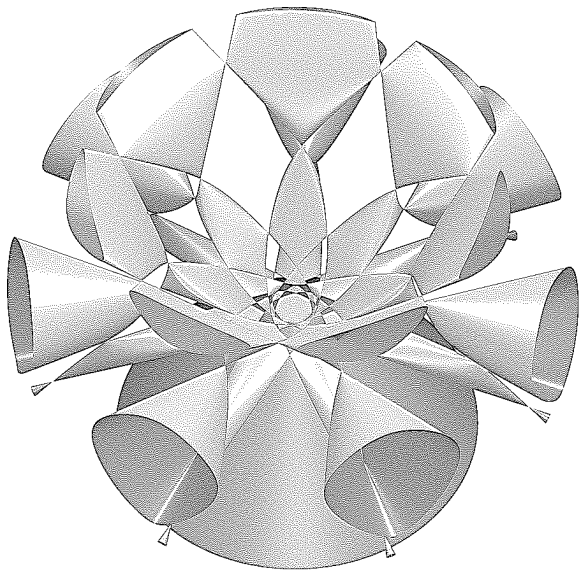}
\caption{A part of the affine chart $w=1$ of the real septic with $99$
  nodes, see \cite{labsAlgSurf} for more images and movies.} 
\label{fig99nodes}
\end{center}
\end{figure}

The existence of the real $\alpha_\dR$ allows us to use the program {\sc surf}
\cite{surf} to compute an image of the $99$-nodal septic $S_{\alpha_\dR}$
(fig.~\vref{fig99nodes}).  
When denoting the maximum number of real singularities a septic in
$\dP^3(\dR)$ can have by $\mu^\dR(7)$, we get, with the remarks mentioned in
the introduction:   
\begin{Cor}
$$99 \le \mu^\dR(7) \le \mu(7) \le 104.$$
\end{Cor}
Note that the previously known lower bounds were reached by S.~V.~Chmutov
($93$ complex nodes: \cite{chmuP3}) and D.~van~Straten ($84$ real nodes: a
variant of Chmutov's construction using regular polygons instead of folding
polynomials).  

As it can be computed using deformation theory and {\sc Singular} that the
space of obstructions for globalizing all local deformations is zero --- this
is based on ideas of D.~van~Straten, details will be 
published elsewhere --- we obtain:  
\begin{Cor}
There exist surfaces of degree $7$ in $\dP^3(\dR)$ with exactly $k$ real nodes
and no other singularities for $k=0, 1, 2, \dots, 99$.
\end{Cor}

Recently there has been some interest in surfaces that do exist over some
finite fields, but which are not liftable to characteristic $0$. 
The reduction of our $99$-nodal septic $S_\alpha$ modulo $5$ (note:
$1\in\dF_5$ satisfies (\ref{eqnalpha}): $7\smcdot1^3 + 7\smcdot1 +1\equiv0$
modulo $5$) neither gives a $99$-nodal surface nor a highly degenerated one as
one might expect because the exponent $5$ appears several times in the
defining equation. 
Instead, we can easily verify the following using computer algebra:
\begin{Cor}
For $\alpha_5:=1\in\dF_5$ the surface $S_{\alpha_5} \subset \dP^3(\dF_5)$ 
defined as in the above theorem has $100$ nodes and no other singularities.
\end{Cor}
Of course, not all the coordinates of its singularities are in $\dF_5$, but in
some algebraic extension. 
The septic has similar geometric properties as our $99$-nodal surface; 
in addition it has one node at the intersection of the $\{x=y=0\}$ axes and
$\{w=0\}$. 
Until now, we were not able to determine if this $100$-nodal septic defined
over $\dF_5$ can be lifted to characteristic zero.

We hope to be able to apply our technique for finding surfaces with many
nodes within families of surfaces to similar problems.
E.g., it should be possible to study surfaces with dihedral symmetry of degree
$9$ and $11$ with many ordinary double points using the same ideas.  
Another application could be the search for surfaces with many cusps.
We already studied families of such surfaces succesfully using computer
algebra in simpler cases \cite{labsSext35}.  

\appendix

\section*{Appendix: The Singular Code to prove the Theorem}
{\small
\begin{verbatim}proc isTrue(int c) { if(c==0) { return("FALSE"); } else { return("TRUE"); } }

ring r = (0,alpha), (x,y,z,w), dp; minpoly = 7*alpha^3 + 7*alpha + 1;
number a(1) = -12/7*alpha^2 - 384/49*alpha - 8/7;
number a(2) = -32/7*alpha^2  + 24/49*alpha - 4;
number a(3) =    -4*alpha^2  + 24/49*alpha - 4;
number a(4) =  -8/7*alpha^2   + 8/49*alpha - 8/7;
number a(5) =    49*alpha^2      - 7*alpha + 50;
poly P = x*(x^6-3*7*x^4*y^2+5*7*x^2*y^4-7*y^6) 
         + 7*z*((x^2+y^2)^3-2^3*z^2*(x^2+y^2)^2+2^4*z^4*(x^2+y^2)) - 2^6*z^7;
poly U = (z+a(5)*w)*(a(1)*z^3+a(2)*z^2*w+a(3)*z*w^2+a(4)*w^3+(z+w)*(x^2+y^2))^2;
poly S = P-U;

"1. Check that the point (1:i:0:0) is not on S:";
poly Si = subst(subst(subst(S, x,1), z,0), w,0); 
ideal yi = y^2+1; yi = std(yi); 
isTrue(reduce(Si, yi) != 0), ": S(1,i,0,0) =", reduce(Si, yi); "";
"2. Check that there is exactly one node on the x=y=0 axes:";
ideal jSxy = x, y, jacob(S); jSxy = std(jSxy);
isTrue(mult(jSxy)==1), ": milnor =", mult(jSxy), ", dim =", (dim(jSxy)-1); "";
"3. Check that milnor(S_y) = 15 (takes some minutes):";
ideal jSy = y, jacob(S); jSy = std(jSy);
isTrue(mult(jSy)==15), ": milnor =", mult(jSy), ", dim =", (dim(jSy)-1); "";
"4. Switch to the affine chart w=1."; "";
S = subst(S,w,1);
ring raff = (0,alpha), (x,y,z), dp; minpoly = 7*alpha^3 + 7*alpha + 1;
poly S = imap(r,S);
"5. Check that all the sing. are in the affine chart w=1 (takes some minutes):";
ideal jSyaff = y, S, jacob(S); jSyaff = std(jSyaff);
isTrue(mult(jSyaff)==15), ": milnor =", mult(jSyaff), ", dim =", dim(jSyaff); "";
"6. Check that all the singularities are nodes (takes approx. an hour):";
poly hessian = det(jacob(jacob(S)));
ideal nonnodes = y, S, jacob(S), hessian; nonnodes = std(nonnodes);
isTrue(mult(nonnodes)==0), ": milnor =", mult(nonnodes), ", dim =", dim(nonnodes); "";}
\end{verbatim}

\nocite{gpSingComAlg}
\nocite{labsAlgSurf}

\bibliographystyle{amsplain} 
\bibliography{papers}

\end{document}